\documentclass[12pt]{article}

\usepackage{amsmath, amssymb}
\usepackage{graphicx}
\usepackage{color}

\usepackage{hyperref}

\hypersetup{
colorlinks=true, 
breaklinks=true, 
urlcolor= blue, 
linkcolor= blue, 
bookmarksopen=true, 
pdftitle={}, 
pdfauthor={Christian Bourdarias}, 
pdfsubject={} 
}

\usepackage{float}

\newcommand{\R}{\mathbb{R}}
\newcommand{\Z}{\mathbb{Z}}
\newcommand{\N}{\mathbb{N}}

\def\1{1\hspace{-.7ex} \rm{I}}

\def\abs#1{\left\vert #1 \right\vert}

\newcommand{\ds}{\displaystyle}

\newcommand{\eps}{ \varepsilon}

\newcommand{\Scal}{\mathcal{S}}

\newcommand{\pr}{{\bf \textit{Proof: }}}
\newcommand{\cqfd}{{\nobreak\hfil\penalty50\hskip2em\hbox{}\nobreak\hfil
$\square$\qquad\parfillskip=0pt\finalhyphendemerits=0\par\medskip}}

\newcommand{\alali}{$\mbox{ }$\\}

\newtheorem{theorem}{Theorem}[section]
\newtheorem{corollary}{Corollary}[section]
\newtheorem{proposition}{Proposition}[section]
\newtheorem{lemma}{Lemma}[section]
\newtheorem{definition}{Definition}[section]
\newtheorem{remark}{Remark}[section]
\newtheorem{ex}{Example}[section]

\title{ \bf Fractional $BV$ spaces 
\\  and  first applications
\\ to scalar conservation laws
}

\author {C. Bourdarias
\thanks{Universit\'{e} de Savoie, LAMA, UMR CNRS 5127, 73376 Le Bourget-du-Lac,
bourdarias@univ-savoie.fr},
M. Gisclon
\thanks{Universit\'{e} de Savoie, LAMA, UMR CNRS 5127, 73376 Le
Bourget-du-Lac, gisclon@univ-savoie.fr}
and S. Junca
\thanks{ Universit\'{e} de Nice Sophia Antipolis, Labo. JAD, UMR CNRS 7351,  Nice,
 junca@unice.fr}
\thanks{Team COFFE, INRIA Sohpia-Antipolis M\'edit\'erann\'ee, 2004 route des lucioles -BP 93, 06902 Sophia-Antipolis, France}
}

\date{\today}

\begin{document}

\bibliographystyle{plain}

\maketitle
\tableofcontents

\abstract{
The aim of this paper is to obtain new fine properties of entropy solutions of nonlinear scalar conservation laws.
For this purpose, we study some  ``fractional $BV$ spaces'' denoted  $BV^s$, for $0 < s \leq  1$,
introduced by Love and Young in 1937.  The   $BV^s(\R)$ spaces are very closed to the critical
Sobolev space  $W^{s,1/s}(\R)$.  We investigate these spaces in relation with  one-dimensional
scalar conservation laws.
$BV^s$ spaces  allow to work with less regular functions than BV functions and appear to be more
natural in this context.  We obtain  a stability result for entropy solutions with $BV^s$ initial
data. Furthermore, for the first time we get the maximal  $W^{s,p}$ smoothing effect  
conjectured by P.-L. Lions,  B. Perthame and E.  Tadmor 
for all  nonlinear degenerate convex fluxes. 
}

 \medskip
  \noindent {\bf AMS Classification}: 35L65, 35L67, 35Q35.

\medskip
 \noindent {\bf Key words}:
 generalized bounded variations, nonlinear  convex flux,
conservation laws, hyperbolicity, entropy solution, Riemann problem.
 {\small}
 
\section{Introduction} 


The space of functions with bounded variation $BV$ plays a key role for scalar conservations laws.
In particular, Oleinik \cite{O57} and Lax \cite{La57} obtained a $BV$ smoothing effect for uniformly
convex fluxes: $f'' \geq \delta > 0$.  
 
Fractional $BV$ spaces, denoted  here  $BV^s$, $0 < s \leq 1 $,  were defined for all $s \in ]0,1[$
in \cite{LY,MO57,MO59,M83}. For $s=1$, $BV^{1}$ is  the space $BV$ of functions
with bounded variation and the space $BV^{1/2}$ is known since $1924$ (\cite{W24}).

Notice that $BV^s$ is not an interpolated space between $L^1$ and $BV$.  Indeed the interpolation between  $L^1$ and $BV$ simply yields $W^{s,1}$ \cite{Tartar}.

The spaces $BV^s$ share some properties with $BV$ and allow to work with less regular functions. 
For the one-dimensional scalar conservation laws,  initial data in $BV^s$ yield  
weak entropy solutions which are still in $BV^s$.  Furthermore, for a degenerate nonlinear convex
flux with only $L^\infty$ data, we obtain  a natural smoothing effect in $BV^s$. Such a smoothing
effect is well known  in the framework of  Sobolev spaces  (\cite{LPT94}). The best
parameter  $s$ quantifying the smoothing effect  is not known in the multidimensional case.  
It is improved in \cite{TT07} and bounded in \cite{DW,Ju11}.  
For the one dimensional case, the best smoothing effect in $W^{s,1}$  conjectured in \cite{LPT94} was first proved in \cite{Ja09}.
We will improve this result in $W^{s,p}$ with $p=1/s$.

It is also well known that the solutions are not $BV$ in the case of a degenerate nonlinear flux,
but they keep some properties of $BV$ functions (\cite{ DOW03,DR}). 
 $BV^s$ spaces appear to be  natural in this context:
 \begin{itemize} 
  \item  we find the maximal $W^{s,p}$  smoothing effect for  a nonlinear degenerate convex flux in one
dimension,
      In this context, $BV^s$ is naturally related to a new one sided H\"{o}lder condition,  
 \item   $BV^s$ spaces share some properties with $BV$ and highlight the $BV$ like structure of
entropy solutions (\cite{DOW03,DR}),
\item  $BV^s$  total variation  is not increasing for all entropy solutions and all fluxes.
 \end{itemize}

In section \ref{sBVs}, we introduce the $BV^s$ spaces and give some usefull properties. We also
investigate for the first time the relations with others classical functional spaces.
In sections \ref{sodscl}  and  \ref{spsa}, we give some applications to scalar conservation laws:
a stability result, the best  smoothing effect in the case of $L^\infty$ data with a degenerate
convex flux  and new results about the asymptotic behavior of entropy solutions  for large time.

\section{ $ BV^s $ spaces} \label{sBVs}

\subsection{Definition}

Let $I$ be a non empty interval of $\R$ and $s \in ]0,1]$. We begin by defining the space $BV^s(I)$ 
which appears to be a generalization of $BV(I)$, space of functions with a bounded variation on
$I$.\\
In the sequel, we note $\Scal(I)$ the set of the subdivisions of $I$, that is the set of finite
subsets $\sigma=\{x_0,x_1,\cdots,x_n\} \subset I$ with   $x_0 < x_1 < \cdots < x_n$.

\begin{definition}
Let be $\sigma=\{x_0,x_1,\cdots,x_n\}\in\Scal(I)$ and let $u$ be a real function on $I$. The
$s$-total variation of $u$ with respect to $\sigma$ is
\begin{eqnarray}
    TV^s u \{\sigma\}  & = &\sum_{i=1}^n |u(x_{i}) - u(x_{i-1})|^{1/s}
\end{eqnarray}
and the $s$-total variation of $u(.)$ on $I$ is defined by

\begin{eqnarray}\label{TVs} 
  TV^s u\{I\}   = \sup_{\sigma \in \Scal(I)} TV^s u\{\sigma\},
\end{eqnarray}

where the supremum is taken over  all the subdivisions $\sigma$ of $I$.
\\
The set $BV^s(I)$ is the set of functions $u:\, I\to\R$ such that $TV^s u\{I\}  < + \infty$.
We define the $BV^s$ semi-norm   by: 
 \begin{eqnarray}\label{snBVs} |u|_{BV^s(I)} =\left( TV^s u\{I\} \right)^s.
 \end{eqnarray}
\end{definition}

We will make use of the following elementary properties:

\begin{proposition}\label{propelem}
Let $I$ be a non empty interval of $\R$  and let $u$ be a real function on $I$.
\begin{enumerate}
 \item For any subinterval $J\subset I$,  $TV^s u\{J\}\leq TV^s u\{I\}$.
\item For any $(a,b,c)\in I^3$ with $a<b<c$, 
$$TV^s u\{]a,b[\}+TV^s u\{]b,c[\}\leq TV^s u\{]a,c[\}.$$
\end{enumerate}
\end{proposition}

\begin{remark}\rm
In the following section it is shown that if $u \in BV^s$ then this function have a finite limit
on the right and on the left everywhere 
(Theorem   \ref{1959}), thus $u$ is measurable and the preceding definition  can be extended
to the class of  measurable functions defined almost everywhere by setting:
$$TV^s u\{I\}~=~\inf_{v=u \, a.e.}  TV^s v\{I\}.$$
\end{remark}

\begin{remark}\rm
For $s=1$, we recover the classical  space $BV(I,\R)= BV^1(I)$.
\end{remark}

\subsection{How to choose a convenient subdivision ?}

In the sequel, we will have to compute explicitly the $s$-total variation of some functions,
especially piecewise constant functions. To this purpose we must know how to get the supremum in
(\ref{TVs}). The following examples and  lemmas show that this calculation can not be done like that of
the total variation in  $ BV $.

\begin{ex}[an increasing function] \label{u=x}
\alali
 Let be $I=[0,1]$, $u(x)=x$ on $I$.
\\
 Then  $ TV^s u([0,1]) = 1$  but   with  the   subdivision 
 $ \ds \sigma_n =  \left \{0, \frac{1}{n}, \cdots, \frac{n-1}{n},1   \right \} $  we have  for all $s<1$, 
$ \ds \lim_{n \rightarrow + \infty} TV^s u\{\sigma_n \} = 0$. 
\end{ex}
So the classical result in $BV$ for smooth function: 
\begin{center}
if $u \in C^1([0,1],\R)$ then  $ TV^s u\{[0,1]\}= \ds \lim_{n \rightarrow + \infty} TV^s u \left \{
\sigma_{n}   \right \} $                          \end{center}
is never  true for all  $s <1$  and for  non-constant function   since
the limit is always $0$.  More generally, refining a  subdivision is not always a good way to
compute the $BV^s$ variation. 

The following example shows two functions with the same $BV$ total variation but never the  same
$BV^s$ total variation for all $s<1$. 
\begin{ex}[a non monotonic function] 
\alali   Let   $a,b$ be some positive numbers,  let  $u$ and $v$ be two functions defined by
$$u=a\,\1_{[0,1[}+(a+b)\,\1_{[1,+\infty[},\qquad v=a\,\1_{[0,1[}+(a-b)\,\1_{[1,+\infty[},$$
where we denote $\1_I$   the indicator function of a set  $I$,  then $TV^s u\{\R\} > TV^s v\{\R\}$ 
for all  $s < 1$.
\end{ex}

This simple phenomenon is related to the monotonicity of $u$ instead of $v$.  We
define  two subdivisions $\sigma_1=\{ -1,\,0,\,1\}$ and $\sigma_2=\{ -1,\,1\}$. We get
easily:
$$TV^s u\{\R\}=  TV^su\{\sigma_2\}=(a+b)^{1/s}>a^{1/s}+ b^{1/s}=  TV^su\{\sigma_1\},$$
$$TV^s v\{\R\}=  TV^sv\{\sigma_1\}=a^{1/s}+ b^{1/s}>\abs{a-b}^{1/s}=  TV^sv\{\sigma_2\},$$
while  $TVu\{\R\}=TV v\{\R\}=a+b= TV u\{\sigma_1\}=TV v\{\sigma_1\}$.
This is an easy consequence of the following lemma, consequence of  the strict convexity of the
function $x \mapsto x^{1/s}$:

\begin{lemma}\label{fondamental}
For all $a,\, b$ in $\R_+^*$ and all $s\in]0,1[$ we have:
 $$\abs{a-b}^{1/s}<  a^{1/s}+ b^{1/s}< (a+b)^{1/s}. $$
 More generally, if $(a_i)_{1\leq i\leq n}$ is a finite sequence of positive real numbers:
 $$ \sum_{1\leq i\leq n} a_i^{1/s} <  \left( \sum_{1\leq i\leq n} a_i \right)^{1/s}.$$
\end{lemma}

To formalize this,  we  propose the following definition:

\begin{definition}
Let  $\sigma = \{ x_0< x_1 <\cdots< x_n \}$ be a subdivision of an interval $I$.  The extremal
points of $\sigma$ with respect to a function $u\, :\, I \rightarrow \R$ are  $x_0$, $x_n$ and, for
$1\leq i\leq n-1$, the points  $x_i$ such that
$\max(u(x_{i-1}),u(x_{i+1} ))  \leq u(x_i)$ or $u(x_i) \leq \min(u(x_{i-1}),u(x_{i+1}))$.
We note $\sigma[u]$ the subdivision of $I$ associated to these extremal points.\\
A subdivision is said to be  extremal with respect to $u$ if $\sigma[u]=\sigma$.
\end{definition}

With this definition, we have the following properties

\begin{proposition}[$BV^s$ variation with extremal subdivisions]\label{sub}~
\begin{enumerate}
\item    
For any subdivision $\sigma$, the  $s$-total variation  of a function $u$ is less or equal to
 $s$-total variation  on the extremal subdivision $\sigma[u]$: 
\begin{eqnarray}
 TV^su\{\sigma \}  & \leq  & TV^su\{\sigma[u]\}, \qquad  \forall \sigma.
\end{eqnarray}
\item
Denote by $Ext(I,u)$ the set of the  subdivisions of an interval $I$, extremal with respect to a
function $u\,:\, I\to\R$. We have 
\begin{equation}\label{Ext} 
  TV^s u\{I\}   = \sup_{\sigma\in Ext(I,u)} TV^s u\{\sigma\}.
\end{equation}
\item 
If u is a monotonic function on the interval $I$ then  $$ TV^s u\{I\} = \left(  \sup_I u - \inf_I u
\right)^{1/s} \hbox{ and } |u|_{BV^s(I)}= TV u \{I\}.$$
\end{enumerate}
\end{proposition}

\pr
\begin{enumerate}
\item
Let $\sigma=\{x_0<x_1<\cdots<x_n\}$ be a subdivision of $I$ and $\sigma[u]=\{
y_0,\cdots,y_N\}$ the subdivision of $I$ associated to the extremal points with respect to $u$. We
introduce the function $\phi : \{0,\cdots,N\} \to\{0,\cdots,n\}$, strictly increasing, such that
 $\phi(0)=0, \, \phi(N)=n$ and $y_j=x_{\phi(j)}$. Setting $u_i=u(x_i)$ we have: 
$$TV^s u\{ \sigma \}=\sum_{i=1}^n \mid u_i-u_{i-1} \mid^{1/s}
=\sum_{j=1}^N  \, \sum_{\phi(j-1) < i \leq \phi(j)}   \mid u_i-u_{i-1}\mid^{1/s}.$$

The sequence $(u_i)$ is monotonic on $[y_{\phi(j-1)},y_{\phi(j)}]$ thus,  by Lemma \ref{fondamental}
$$\sum_{\phi(j-1) < i \leq \phi(j)}   \mid u_i-u_{i-1} \mid^{1/s} \leq  \mid
u_{\phi(j)}-u_{\phi(j-1)} \mid^{1/s}.$$

Finally, $\ds TV^s u\{ \sigma \} \leq \sum_{j=1}^N 
\mid u_{\phi(j)}-u_{\phi(j-1)}\mid ^{1/s} =TV^s u \{\sigma[u]\}$.
\item It is  a direct consequence of the first item of Proposition \ref{sub}. 
\item   The extremal subdivision  for a monotonic function have only two points: $ \sigma =\{\min
\sigma, \max  \sigma\}$ 
 and the result follows. 
\end{enumerate}
\cqfd

We have seen in Example \ref{u=x} that we can have $\tau \subset \sigma$ but $TV^s u\{\tau\} >TV^s
u\{\sigma \} $: take $\tau=\{0,1\}$ and $\sigma=\sigma_n$ with $n >1$.
The following example shows that this problem can also occur for extremal subdivisions.

\begin{ex}[A piecewise monotonic function]\label{cex}
Let  $I=[0,3]$,   let  $w$ be the  continuous piecewise linear function defined by: 
 $w(0)=0$, $w(1)=a$, $w(2)=a- \eps$,  $w(3)=b$,  with $0 < \eps < a < b$, let
  $\tau = \{0,3\}$ and  $\sigma =\{0,1,2,3\}$.  
$\tau$ and $\sigma$  are extremal subdivisions, we have $\tau \subset \sigma$ but $TV^s u\{\tau\}
>TV^s u\{\sigma \}  $ for  all $s<1$ and $0 < \eps $ small enough.
\end{ex}

Indeed we have    $$TV^s u\{\tau\}=b^{1/s} >  a^{1/s}+ (b-a)^{1/s}=TV^s u\{\sigma \}  $$ and 
$$TV^su\{\sigma \} =  a^{1/s}+  \eps^{1/s}+ (b - a + \eps)^{1/s} = g(\eps) .$$ 
We have also  $g(0) =  a^{1/s}+ (b-a)^{1/s}< TV^s u\{\tau\}=b^{1/s}$ 
by Lemma \ref{fondamental} and $g$ is a continuous function, thus
the inequality holds for $0< \eps$ small enough. 
\\
Example \ref{cex} shows that the $TV^s$ variation of a function is not necessarily computed 
 using all extremal points of this function.
 \\
 Conversely, the following proposition is useful to compute  $BV^s$ variation of oscillating
functions  with diminishing amplitudes. 
  
 \begin{proposition}[ $BV^s$ variation of alternating diminishing oscillations]\label{ADO}~\\
Let $I= \bigcup_{k \geq 0} I_k$, $I_k=[x_{k},x_{k+1}[$, $x_{k} <x_{k+1}$ and $u$ be 
 a monotonic func\-tion on each $I_k$,  with  successive different monotonicity: 
 $(u(x) - u(y)) (u(z) - u(t)) \leq 0$ for all $ x_k \leq x < y < x_{k+1} \leq z < t < x_{k+2}$.
The oscillation of  $u$ on the compact interval $\overline{I_k}$ is
 $a_k = \sup_{ x,y \in [x_k,x_{k+1}]} |u(x) - u(y)|$.
 
 If the oscillation  $(a_k)_k$  is monotonic then
 
 \begin{eqnarray}\label{TVsmonotonic}
 TV^su\{I \}  & =   &  \sum_k a_k^{1/s}.
\end{eqnarray}
  \end{proposition}
  
 Notice that if the sequence of successive amplitudes is not monotonic   then \eqref{TVsmonotonic}
can be wrong.  The result is still valid with a  finite union of $I_k$. For  the infinite case, the
non increasing oscillations is the interesting case. In this case, the proposition states:
  
$$ u \in BV^s(I)  \mbox{ if and only if  }(a_n) \in l^p(\mathbb{N}) \mbox{ with  } s\;  p=1.$$ 

\pr
  to prove that   $TV^su\{I \}  =  \sum_k a_k^s$  we restrict ourselves  to the case of a
piecewise constant  function. The general case follows. 
   \\
Let   $A_N = a_0 - a_1+ \cdots + (-1)^N a_N$  and $u(x) = A_N$ on $I_N$. 
The inequality   $TV^su\{ I \}  \geq  \sum_k a_k^s$  is clear by taking the subdivision $\sigma*
=\{x_0, x_1,\cdots \}$. Let   $\sigma= \{y_0, y_1,\cdots \}$ be any  other extremal
subdivision.   We can assume that there is at most one $y_j$ in each $I_n$ since the contribution
is zero for two extremal points in the same interval.
  
  Let us define $k(j)$ by the condition $x_j \in I_{k(j)} $. We have  to prove that 
  $$
    \sum_{j} |u(y_{j+1}) - u(y_{j})|^p \leq  \sum_{i} |u(x_{i+1}) - u(x_{i})|^p.
  $$
  We have $ |u(y_{j+1}) - u(y_{j})| \leq  |u(x_{k(j)+1}) - u(x_{k(j)})|$ since $(A_N)$  is the
partial sum of an alternating series.  This is enough to conclude the proof.
 \cqfd

Let us study a more complex example where the sequence of the increasing jumps belongs to $l^1$ and
the sequence of the decreacreasing jumps belongs to $l^2$.  Does the function  belong to
$BV^{1/2}$? The result is more surprising.

\begin{ex}\label{oscillatingBOOM}
Let $( a_n)_n$ be   a positive sequence  which belongs  to  $l^1(\N)$ such that  $b_n= \sqrt{a_n}$
does not  belong to  $l^1(\N)$. We set:
 $$  z(x) = \sum_{n} z_n \1_{]n-1,n]}(x), \quad    z_{2n+1}= z_{2n} +a_n,\;    z_{2n+2}= z_{2n+1} -
b_n,\; z_0=0. $$
  $z$ does not  belong to any $BV^s(\R)$ for all $s$. 
\end{ex}

\noindent  \underline{Proof 1}: notice that $ \sum a_n < \infty $ and  $ \sum b_n^2  < \infty $
since $b_n^2=a_n \in l^1$.\\
$z_{2n+2} = (a_0 + \cdots + a_n) - (b_0 + \cdots + b_n ) $ yields  $\ds   \lim_{n \rightarrow +
\infty} z_{2n}=- \infty$ and also  $\ds   \lim_{n \rightarrow + \infty} z_{n}=- \infty$ .
 This implies  $\ds \lim_{x \rightarrow + \infty} z(x)= - \infty$:  $z $ is not bounded  and
thus in none $BV^s$ thanks to Proposition \ref{pBVs} below.
\cqfd
 
\noindent \underline{Proof 2}:
notice that  $a_n = o (b_n)$  
 and  for $n$ large enough     $TV^su\{ ]2n+1, 2n +3[\} \sim ( b_n + b_{n+1})^{1/s}$ 
  in a similar way as in  Example \ref{cex}.     
For any $k >0$,  in a same way,  we have  $TV^su\{ ]2n+1, 2n +2 k +1[\} \sim ( b_n + \cdots +
b_{n+k})^{1/s}$, but $\sum_n b_n = + \infty$ so the  $BV^s$ total variation blows up. 
\cqfd 

\noindent \underline{Proof 3}:
 there is  another way to interpret Example \ref{oscillatingBOOM}. Functions $L^\infty$ with
total increasing varition bounded are $BV$.  By construction,  the total increasing variation  $TV_+
z = \sum_n a_n $ is bounded, but $z$ is not $BV$  since 
 the  total decreasing variation is not bounded: $TV_- z = \sum_n b_n = + \infty$.  So, $z$ is not
in $L^\infty$ and also in none $BV^s$.  
 \cqfd 

The problem is more complicated if  we assume that  $(a_n)_n$  does not belongs to $l^1$ . 
The previous  argument in $BV$  is not known in $BV^s$ for $s<1$.  For instance, if  $(a_n)_n$  does
not belong to $l^1$ but belongs to   $l^2$, is $z $ in $BV^{1/4}$ ?

\subsection{Some properties of $BV^s$ spaces}\label{properties1}

We begin with some properties of $BV^s(I)$ which arises directly from the definition:

\begin{proposition}\label{pBVs}
Let $I$ be an  interval  of $\R$. The following inclusions hold:
\begin{enumerate}
\item for all $s\in]0,1]$,  $BV^s(I) \subset  L^{\infty}(I)$,
\item  if $0 < s < t \leq 1$  and $I$ is not reduced to one point then $BV^t(I) \subsetneqq
BV^s(I)$.
\end{enumerate}
\end{proposition}
\pr
\begin{enumerate}
\item  
Let $a \in I$. For any $x \in I$ one has $|u(x) - u(a)| \leq |u|_{BV^s(I)} $ thus 
$ \|u\|_{L^{\infty}(I)}\leq  |u(a)|+ |u|_{BV^s(I)}$
then the first inclusion holds.
\item 
 We can  assume $I=]0,1[$ without loss of generality. 
The null function of course belongs to all spaces $BV^s$. Assume $u\neq 0$ and  $u$ in
$BV^t(I)$ for some $t\in]0,1]$ and let $s$ be such that $0 < s < t$.
First, $u \in L^\infty(I)$ and $v=\ds \frac{u}{2\|u\|_\infty}\in BV^t$. 
Now $\|v\| _\infty = 1/2$ thus  for any variation $\Delta v$ of $v$ we have $|\Delta v|\leq 1$ and
$|\Delta v|^{1/s} \leq |\Delta v|^{1/t}$ then the second inclusion follows.
 
In order to prove that $BV^t(I) \neq BV^s(I) $, let us consider the  function 
$$u(x) =\sum_{n=1}^{+\infty} a_n \1_n(x)$$
where  $\1_n$ is the indicator function of $I_n=](n+1)^{-1},n^{-1}]$ and $a_n=\ds\sum_{p=1}^n
\frac{(-1)^p}{p^t}$.\\
On one hand, choosing the subdivision $\ds \sigma_n=\{\frac{1}{p}\, ;\, 1\leq p\leq n\}$ (extremal
with respect to $u$) we get $TV^tu\{]0,1[\}\geq \ds\sum_{p=1}^n \frac{1}{p}$ and  $u\notin
BV^t(]0,1[)$. On the other hand, using the same family of subdivisions $\sigma_n,\,n\geq 1$ and
Proposition \ref{ADO} we get, for $0<s=t-\eps$,  $ TV^{s}u\{]0,1[\}=\ds\sum_{n=1}^\infty
\frac{1}{n^{\frac{t}{t-\eps}}}<+\infty$ thus $u\in BV^{s}(]0,1[)$.
\end{enumerate}
\cqfd

\begin{proposition}\label{lim+-}
If $u \in BV^s(I)$  then $u$ is a regulated function.
\end{proposition}
This result is already found in \cite{MO59}.  We give a proof for the convenience of the reader.
\\
\pr 
let be $(a, \, b) \in I^2$ with $a<b$, let $\eps>0$ and  $\sigma\in\Scal(]a,b[)$ be such that $$
TV^s u \{\sigma\}\geq  TV^s u \{]a,b[\} -\eps.$$
There exists $\alpha>0$ such that $\sigma\in\Scal(]a+h,b[)$ for any $h\leq\alpha$ and we have 
(Proposition \ref{propelem}):
$$TV^s u \{]a,a+h[\}+TV^s u \{]a+h,b[\}\leq TV^s u \{]a,b[\}, $$
thus $h\leq\alpha$ implies  $ TV^s u \{]a,a+h[\}\leq\eps$ i.e. $\ds\lim_{h\to 0}TV^s u
\{]a,a+h[\}=0$.
The oscillation of $u(.)$ on $]a,a+h[$ also tends to $0$ as $h\to 0 $ and this is enough to get a
right limit for $u$ at point $a$.
 For the existence of a left limit and  the cases of  $\alpha=\inf I$ and  $\beta=\sup
I$, the proof is very similar.\cqfd

\bigskip

Recall that for $\alpha>0$ and $p\geq 1$ a function $u$ belongs to the space
$Lip(\alpha,L^p(\R))$ if there exists some constant $c\geq0$ such
that $||u(\cdot+h)-u||_{L^p} \leq c \, |h|^\alpha$ for all $h\in\R$ (\cite{De98}). The space
$BV(\R)$ is nothing but $Lip(1,L^1(\R))$ and if $u\in BV(\R)$ we have 
$$TV(u)=\sup_{h > 0}\frac{1}{h} \int_\R |u(x+h)-u(x)|\, dx.$$

Dealing with the space $BV^s(\R)$, we have a  different result:

\begin{proposition} \label{LipLp}
For any $0<s<1$,  $BV^s(\R) \subset  Lip(s,L^{1/s}(\R))$. More precisely,
for $u\in BV^s(\R)$: 
\begin{equation}\label{ineqLip}
\sup_{h > 0} \frac{1}{h} \int_\R |u(x+h)-u(x)|^{1/s} dx 
 \leq    \ds TV^su\{\R\}
\end{equation}
and this inequality generally cannot be replaced by an equality.
\end{proposition}

\pr
for $u \in BV^s(\R)$ and $h>0$ we have:
\begin{eqnarray*} 
\ds \int_\R |u(x+h)-u(x)|^{1/s}  \, dx  & = & 
\sum_{k \in \Z}  \int_{kh}^{kh +h} |u(x+h)-u(x)|^{1/s} \, dx \\
& = & 
\sum_{k \in \Z}  \int_{0}^{h} |u((k+1)h+y)-u(kh+y)|^{1/s} \, dy \\
& = & 
\int_{0}^{h}  \sum_{k \in \Z}   |u((k+1)h+y)-u(kh+y)|^{1/s} \, dy \\
&  \leq &  \int_{0}^{h} TV^su\{\R\}  \, dy   =  h  \, TV^su\{\R\}.
\end{eqnarray*}
\medskip

Inequality (\ref{ineqLip}) and the inclusion $BV^s(\R) \subset  Lip(s,L^{1/s}(\R)))$ follow.

In order to prove that Inequality (\ref{ineqLip}) may be strict, we consider  the
function $u(x)= x \1_{[0,1]}$ and we set, for $p\geq 1$ and $h>0$: 
$$I_p(h)=\ds  \frac{1}{h}\int_\R |u(x+h)-u(x)|^{p} \, dx.$$
On one hand $TV^su\{\R\}=2$. On the other hand:

if $h \geq 1$,  then $$h \, I_p(h)=\int_{-h}^{1-h}  |x+h|^{p} \,  dx + \int_{1-h}^{0}  0  \, dx +
\int_{0}^1 x^p \, dx= \frac{2}{p+1},$$
thus  $I_p(h)\ds \leq\ I_p(1)= \frac{2}{p+1}<2 $,

if $0<h \leq 1$,  then 
\begin{eqnarray*} 
h \, I_p(h)&=&\int_{-h}^0   |x+h|^{p} \,  dx + \int_0^{1-h}  h^p \, dx + \int_{1-h}^1 x^p \, dx\\
&=&\frac{ h^{p+1}}{p+1} + h^p (1-h) + \frac{1- (1-h)^{p+1}}{p+1},
\end{eqnarray*}
and in particular $I_1(h)=2-h$. For $p>1$  $I_p(0^+)=1$, thus  $\ds\sup_{h>0}I_p(h)=1$ or there
exists $h_0> 0$ such that $\ds\sup_{h>0} I_p(h)=I_p(h_0)$. Now,  $I_p(h_0)$ is non increasing with
respect to $p$ because $ |u(x+h)-u(x)|\leq 1$ thus $I_p(h_0)\leq I_1(h_0)<2$. Finally we get:
$$ \sup_{h > 0} \frac{1}{h} \int_\R |u(x+h)-u(x)|^{1/s} dx = \sup_{h > 0}  I_{1/s} (h) < 2 = TV^s
u\{\R\}.$$
\cqfd

\begin{corollary}
For any $0<s<1$ and any interval $I\subset\R$ (with $\mathring{I}\neq \emptyset$) we have
$$BV^s(I) \subset Lip(s,L^{1/s}(I)).$$
Moreover, with $I_h = \{x \in I, \; \mbox{ such that }  x+h \in I\}$, we have:
$$ \sup_{h>0} \frac{1}{h} \int_{I_h} |u(x+h)-u(x)|^{1/s} dx \leq TV^su\{I\},$$
and this inequality generally cannot be replaced by an equality.
\end{corollary}

\pr this result follows immediately from Proposition \ref{LipLp} thanks  to the following
lemma.\cqfd

\begin{lemma}
 Let $I\subset\R$ be  an interval. We set $a=\inf I$ and $b=\sup I$. For $u\,:\, I\to\R$ we note
$\tilde{u}\,:\, \R\to\R$ the extension of $u$ such that:
\begin{description}
 \item[-] if $a\in I$ then $\tilde{u}(x)=u(a)$ for $x\leq a$,
  \item[-] if $a\notin I$ and $a\neq -\infty$ then $\tilde{u}(x)=u(a^+)$ for $x\leq a$,
  \item[-] if $b\in I$ then $\tilde{u}(x)=u(b)$ for $x\geq b$,
  \item[-] if $b\notin I$ and $b\neq +\infty$ then $\tilde{u}(x)=u(b^-)$ for $x\geq a$, 
 \end{description}
 then  
 $$ \sup_{|h|<dist(x,\partial I)} \frac{1}{h} \int_I |u(x+h)-u(x)|^{1/s} dx \leq \sup_{h\neq  0}
\frac{1}{h} \int_\R |\tilde{u}(x+h)-\tilde{u}(x)|^{1/s} dx$$ 
and $TV^s\tilde{u}\{\R\}= TV^s u\{I\}$.
\end{lemma}

\pr the first inequality is obvious. Next, on one hand we have trivially $TV^s u\{I\}\leq TV^s
\tilde{u}\{\R\}$. On the other hand, in
order to get the converse inequality it suffices to consider the case
$I=]-\infty,b]$. Let $\tau\in\Scal(\R)$ be such that $\sigma\tau\cap I\neq\emptyset$ and $\tau\cap
I^c\neq\emptyset$. If $\sigma=\{x_0<\cdots<x_n\}$ then we get easily:
\begin{eqnarray*}
 TV^s\tilde{u}\{\tau\} & = & TV^s\tilde{u}\{\sigma\cup\{x_{n+1}\}\}\\
 & = & TV^s\tilde{u}\{\sigma\cup\{b\}\}\\
 & = & TV^s u \{\sigma\cup\{b\}\} \leq  TV^s u \{I\},
\end{eqnarray*}
thus $TV^s\tilde{u}\{\R\}\leq TV^s u\{I\} $.
\cqfd

\bigskip

Some results of the next proposition can be found in \cite{MO59}.
There are the same properties for the space $BV$.

\begin{proposition} \label{1959}
Space  $BV^s(I)$ is endowed with the  following properties:
\begin{enumerate}
\item 
$BV^s(I) \cap L^{1/s}(I)$ with the norm $\|u\|_s = \|u\|_{L^{1/s}}  +
|u|_{BV^s(I)}$ is a Banach space,
\item 
the embedding $BV^s(I) \cap L^{1/s}(I)\hookrightarrow L^1_{loc}(I)$ is compact. 
\end{enumerate}
\end{proposition}

\pr 
\begin{enumerate}
\item
the proof is classic (\cite{MO59}).
\item 
Case $I=\R$: It suffices to prove that  $BV^s(\R) \cap L^{1/s}(\R)$ is compactly imbedded in
$L^{1/s}_{loc}(\R)$ because  $L^{1/s}_{loc}(\R)\hookrightarrow L^1_{loc}(\R)$.  
This is a direct consequence of the  Riesz-Fr\'echet-Kolmogorov Theorem  since 
$(u_n)$ is bounded in $BV^s$  and we have from Proposition \ref{LipLp}: 
\begin{eqnarray}
\ds  \int_\R |u_n(x+h)-u_n(x)|^{1/s} dx & \leq & |h|  \, TV^su_n\{\R\}  \leq C |h|.
\end{eqnarray}
The proof is similar in the general case (see for example \cite{MO59}). 
\end{enumerate}
\cqfd 

To end this section,we give two approximation results which will be usefull in the context of scalar
conservation laws (see Section \ref{sodscl} below). 

\begin{proposition} \label{approx}
Let $I$ be an interval of $\R$ and let $u$ be a function in $BV^s(I)$. There exists a sequence $(u_n)_{n\geq 0}$ of step functions such that $u_n\to u$ in $L^1_{loc}$ and $TV u_n\{I\}\leq TV u\{I\}$.
\end{proposition}

\pr
we treat the case $I=\R$ for the sake of simplicity. Let $h>0$, we set $u^h=\ds\sum_p u_p^h\,
\1_{]ph,(p+1)h]}$ with $\ds u_p^h=\frac{1}{h}\int_{ph}^{(p+1)h}u(x)\,dx$: we have $u^h\to u$ in
$L^1_{loc}$ as $h\to 0$.\\
For all $p\in\N$ there exists $x_p^h,\, y_p^h \in ]ph,(p+1)h[$ such that $u(x_p^h)\leq u_p^h\leq
u(y_p^h)$. Let us consider a maximal finite sequence $p_i, \, p_i+1,\cdots,p_{i+1}$ corresponding to
a monotonic  sequence $(u_{p_i}^h,\cdots,u_{p_{i+1}}^h)$: we set $x_i=x_{p_i}^h$ if the sequence is
increasing, $x_i=y_{p_i}^h$ else. The maximality of the sequence of indexes ensures the consistency
of this definition. Let $\sigma$ be a subdivision $\{x_j<x_{j+1}<\cdots<x_{j+k}\}$. By Lemma
\ref{fondamental} we have clearly $TVu\{\R\}\geq TVu[\sigma]\geq TVu^h[\sigma]$ and thus
$TVu\{\R\}\geq TVu^h\{\R\}$. Proposition \ref{approx} follows immediately.
\cqfd

\begin{proposition} \label{TVSCI}
Let $(u_n)_{n\geq 0}$ be a sequence of $BV^s(\R)$ functions such that $u_n\to u$ a.e., then $TV^su\{\R\}\leq\liminf TV^su_n\{\R\}$.
\end{proposition}

\pr let $\sigma=\{x_0<x_1<\cdots<x_p\}$ be a subdivision of $\R$. We have
$TV^su_n[\sigma]=\ds\sum_{i=1}^p |u_n(x_{i}) - u_n(x_{i-1})|^{1/s}\to TV^su[\sigma]$ as $n\to\infty$
and $TV^su_n[\sigma]\leq  TV^su_n\{\R\}$. Thus $TV^su[\sigma]\leq \liminf TV^su_n\{\R\}$ and the
result follows.
\cqfd


\subsection{ Relations between $BV^s$ and  $W^{s,1/s}$}\label{properties2}


Fractional Sobolev spaces are used in \cite{LPT94} to study the smoothing effect for nonlinear
conservation laws. An aim of this paper is to show that $BV^s$ space are more appropriate  to study
the smoothing effect for nonlinear conservation laws.

Let us  first compare  $BV^s$ and $W^{s,p}$. 
Roughly speaking  $BV^s \simeq W^{s,1/s}$ but $BV^s\neq W^{s,1/s}$. More precisely  $W^{s,p}$, when 
$ \ds   s \, p = 1$  is the borderline Sobolev space in dimension one.
Indeed the embedding in the space of continuous function just fails:
 \begin{itemize}
  \item    $ p > \frac{1}{s}    \implies W^{s,p}(-1,1) \subset C^0([-1,1])$,
  \item     $ p < \frac{1}{s}    \implies        H \in W^{s,p}(-1,1) $ where $H$ is the
Heaviside step function,
  \item   For $ p = \frac{1}{s}$,  $H \notin W^{s,1/s}(-1,1)$, 
  but  some more complicated discontinuous functions are in 
   $ W^{s,1/s}(-1,1)$  such that $\ln \ln |x|$ which is not bounded and $\sin \ln \ln |x|$  which is
bounded but discontinuous (\cite{BN}).
 \end{itemize} 
For the classical $BV$ space endowed with the norm: $ \|u\|_{BV} = \|u\|_{L^1} + TV u$ we have: 
$$ W^{1,1}(\R) \subsetneqq BV(\R) \subsetneqq   \bigcap_{s <1 } W^{s,1}(\R) .$$ 
\begin{proposition}[$BV^s$ and $W^{s,p}$]\label{pSobolevBVs}
\alali
Let $I\subset\R$ be a  nontrivial bounded interval, then 
\begin{enumerate}
\item $W^{s,\infty}(I) \subset BV^s(I)$,
\item  $ \ds BV^s(I)  \subset    \bigcap_{t <s }  W^{t,1/s}(I)$,
\item $BV^s(I) \neq W^{s,1/s} $,   more precisley  we have 
$BV^s(I) \nsubseteq W^{s,1/s} $,    $BV^s(I) \nsupseteq W^{s,1/s} $.
%
\end{enumerate}
\end{proposition}

\pr  ~
\begin{enumerate}
\item 
Let $u \in W^{s,\infty}(I)$. There exists $C>0$ such that $|u(x)-u(y)| \leq C |x-y|^s$, 
$$ TV^s u \{\sigma\}   =   \ds \sum_{i=1}^n |u(x_{i}) - u(x_{i-1})|^{1/s} 
 \leq  C^{1/s}\sum_{i=1}^n |x_{i} - x_{i-1}| \leq  C^{1/s}  |x_n -x_0| $$ and $u \in BV^s_{loc}$.   

\item
An usual semi-norm on fractional Sobolev space is:  
\begin{equation}\label{Wsp}
 |u |^{p}_{W^{s,p}(\R)}=\int _{\R} \int_{\R} \frac{  |u(x)-u(y)  |^p}{  |x-y|^{1+sp}} \, dx \, dy
= \int_{\R} \int_{\R} \frac{ |u(x+h)-u(x) |^p}{ |h |^{1+sp}} \, dx \, dh.
\end{equation}

Now, assume $u \in BV^s(I)$. We note  $p=1/s$.
We bound $|u|_{\sigma}^p $ the intrinsic semi-norm    of
 $W^{\sigma,p}(I)$ by:
 \begin{eqnarray*}  |u|_{\sigma}^p  & = &
   \ds \int_{-l}^l  \int_{a}^{b-h}  \frac{|u(x+h)-u(x)|^p}{|h|^{p\sigma+1}} dx \, dh  \\
   & \leq & \ds \int_{-l}^l \left( \frac{1}{|h|} \int_a ^{b-h}  |u(x+h)-u(x)|^p dx        \right)
\frac{dh}{|h|^{p\sigma}}  \\
 &  \leq &  \ds TV^su \int_{-l}^l  \frac{dh}{|h|^{p\sigma}} < + \infty
\end{eqnarray*} 
thanks to Poposition \ref{LipLp}
and because $p \,  \sigma = p(s-\varepsilon)=1- \, p  \, \varepsilon <1$.
\item 
 More precisely there is no inclusions between $BV^s$ and $W^{s,1/s}$. 
 \begin{enumerate}
\item 
$W^{s,1/s}$ is not a subspace of  $BV^s$: the Heaviside function is in $BV^s$ but not in $W^{s,1/s}$: use the integral criterium
(\ref{Wsp}).
\item
$BV^s$ is not a subspace of  $W^{s,1/s}$:  we have just to consider the following example (cf
\cite{BN}):

$\ln |\ln|x||  \in W^{s,1/s} \mbox{ but }  \ln |\ln|x|| \notin BV^s \mbox{ (it is not bounded). }$
\end{enumerate}
\end{enumerate}
\cqfd

 \begin{remark}
  $BV^s(I)$ is not a Sobolev space.
 \end{remark} 
 
That is to say,  the set $BV^s$  (resp. the $TV^s$ variation) is not a Sobolev space (resp. its
semi-norm).
 Indeed,  the $s$-total variation is invariant under dilations.
Indeed, for any $\lambda \neq 0$, the function $u_\lambda$ defined by  $u_\lambda(x)=u(\lambda x)$
satisfies  $TV^s u_\lambda\{\R\}=TV^s u\{\R\}$. 
Thus for compactly spported functions the $s$-total variation is independent of the
support. In par\-ti\-cu\-lar, the $s$-total variation is not related to  a Sobolev semi-norm except
for $W^{s,p} $ with $sp=1$.
 But $W^{s,1/s}$ and $BV^s$ are different. So the $s$-total variation is not a Sobolev semi-norm. 
 

\section{$BV^s$ stability for scalar conservation laws} \label{sodscl} 


\begin{theorem}\label{stab}
Let $u_0 \in BV^s(\R)$, $f \in C^1(\R,\R)$  and $u$
be the unique entropy solution  on $]0,+\infty[_t\times\R_x$ of
\begin{eqnarray} \label{eqscl} 
  \partial_t u + \partial_x f(u) = 0, 
 & \quad & u(0,x)=u_0(x), 
\end{eqnarray}
then
\begin{equation}\label{TVSD} 
\forall t>0\quad TV^s u(t,.)(\R) \leq TV^s u_0(.)(\R).
\end{equation}
\end{theorem}

\begin{figure}[H]
\centering
\includegraphics[scale=0.6]{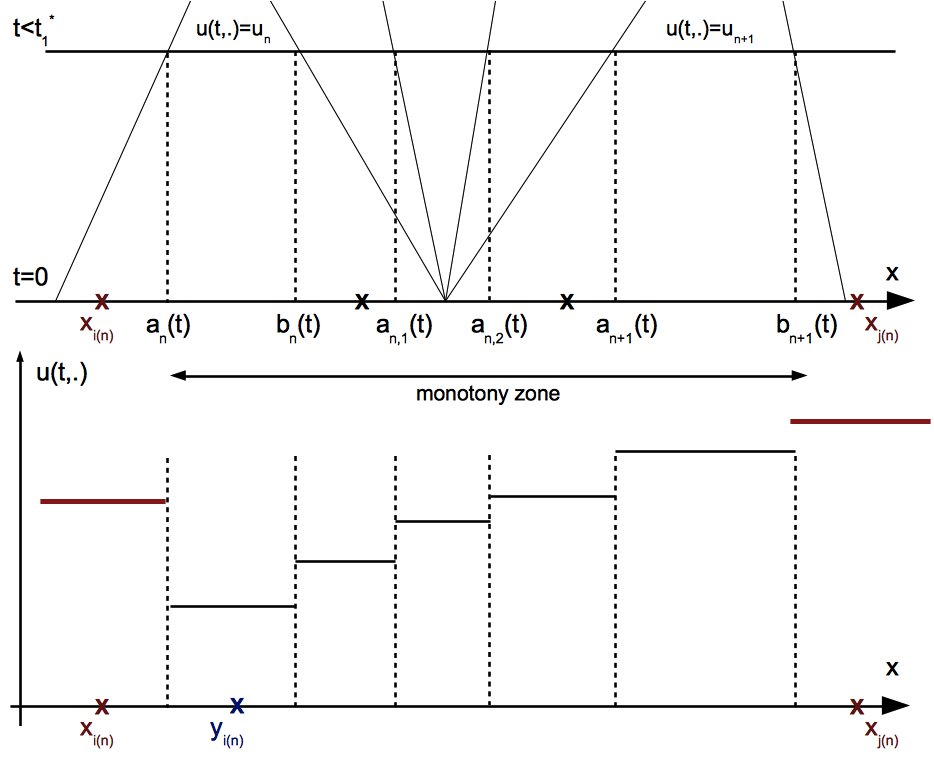} 
\caption{\small\it  a zoom around a wave fan. The $\times$ symbols correspond to a subdivision in
the neighborhood of a wave fan: we have here $P=\{y_{i(n)}\}$.}\label{stefan}
\end{figure}

This theorem means that the $s$-total variation is not increasing with respect to time.
\\[2mm]
\pr 
in a first step we show that this property is achieved for an approximate solution obtained
with the Front Tracking Algorithm (\cite{Br00,D00}), thus we assume  that the initial
condition is piecewise constant and writes $u(0,x)=u^0(x)=\ds\sum_n u_n^0 \1_{]a_n,b_n]}(x)$.
The key point is that  the solution of the Riemann problem at each point of discontinuity,
consisting in a composite wave, is  piecewise constant and  monotonic. Actually, in the framework
of the Front Tracking Algorithm, we also assume that the flux function $f$ is piecewise affine, thus
we have to deal with $K$ contact discontinuities for each Riemann Problem, where $1+K$ is the number
of intervals where $f$ is affine. In a second step we show that we can pass to  the limit in this
approximation process in order to get (\ref{TVSD}).\\[2mm]
\textbf{First step - }
We denote by $t_1^*$ the time of the first interaction and, following \cite{Br00},  we can suppose
that there exists an only interaction. For $t < t_1^*$, we denote 

$$u(t,\cdot)=\sum_n \left(u_n \1_{]a_n(t),b_n(t)]}+
\sum_{m=1}^K u_{n,m} \1_{]a_{n,m}(t),a_{n,m+1}(t)]}\right),$$
 
 where $b_n<a_{n,1}<\cdots<a_{n,K+1}<a_{n+1}$  (zone corresponding to the wave fan denoted $F_n$:
see Fig. \ref{stefan}), with the monotony condition:

$$u_n \leq u_{n,1} \leq \cdots\leq u_{n,K} \leq u_{n+1} \hbox{ or }u_n \geq u_{n,1} \geq \cdots\geq
u_{n,K} \geq u_{n+1}$$

Let $\sigma=\{x_0,\cdots ,x_p\}$ and $TV^su(t,.) \{ \sigma \}=\sum \mid
u(t,x_i)-u(t,x_{i-1})\mid^{1/s}$. Let $\tilde{\sigma}$ be the subdivision obtained by removing  the
points $x_i$ located in a fan zone: $\tilde{\sigma}=\sigma\setminus
\ds\bigcup_n[b_n(t),a_{n+1}(t)]$. We are going to show that it is possible to add to
$\tilde{\sigma}$ a finite set $P$ of points located in $ \ds\bigcup_n ]a_n(t),b_n(t)[$ in such a way
that $TV^su(t,.) \{  \tilde{\sigma}\cup P \} \geq TV^su(t,.) \{ \sigma \}$. This being carried out,
we get $TV^su(t,.) \{ \sigma \}\leq TV^su(t,.) \{ \tilde{\sigma}\cup P \}\leq TV^su_0 $ and thus
(\ref{TVSD}) holds for the \emph{exact solution} of Problem (\ref{eqscl}) associated to the
\emph{approximate initial condition} and the \emph{approximate (piecewise affine) flux}.

\medskip

In the bounded interval $[\min \sigma, \max \sigma]$ there is a finite number of fan zones and we
have just to consider the case of a single wave fan $F_n$ and the associated monotony zone
$M_n=]a_n(t),b_{n+1}(t)[ $ in which we assume (for instance) that $u(t,\cdot)$ is increasing.\\ 
If $\sigma\cap  M_n=\emptyset$ then we have nothing to do, else we set  $i(n)=\max\{0\leq i\leq p\,
;\, x_i\leq b_n(t)\} $ (if exists) and $j(n)=\min\{0\leq i\leq p\, ;\, x_i\geq a_{n+1}(t)\} $ (if
exists).

\begin{itemize}
\item If $i(n)$ exists and $u(x_{i(n)})>u_n$ then we add to $\tilde{\sigma}$ any point  $y_{i(n)}\in
]a_n(t),b_n(t)]$,
\item if $j(n)$ exists and $x_{j(n)}<u_{n+1}$ then we add to $\tilde{\sigma}$ any point\\
$y_{j(n)}\in ]a_{n+1}(t),b_{n+1}(t)]$, else we have nothing to do.
\end{itemize}
Let $P$ be the set of the added points according to the preceding procedure. Thanks to Lemma
\ref{fondamental}, we get immediately $TV^su(t,.) \{ \sigma^* \} \geq TV^su(t,.) \{ \sigma \}$ 
where $\sigma^*=\tilde{ \sigma}\cup P $.\\
When the first interaction occurs ($t=t_1^*$),  it appears a new monotony zone where the solution 
varies between two successive values taken by $u(t,\cdot)$ for $t$ in some interval
$[t_1^*-\epsilon,t_1^*[ $, thus the total variation does not increase.This concludes the first step.

\textbf{Second step - }
Let $(u_0^n)_{n\geq 0}$ be a sequence of step functions in $BV^s$  such that $u_0^n\to u_0$ in
$L^1_{loc}$ and a.e., with $TV^s u_0^n\leq TV^s u_0$ : this is ensured by Proposition
\ref{approx}.\\
Let $(f^n)_{n\geq 0}$  be a sequence of piecewise affine functions such that $f^n\to f$ uniformly on
every compact set.\\
Let $u^n$  be the solution of Problem (\ref{eqscl}) associated to the initial condition $u_0^n$ and
the  flux $f^n$. For all $t\geq 0$, $(u^n(t,\cdot))_{n\geq 0}$ is bounded in $L^\infty\cap BV^s$ 
thus it converges,  extracting a subsequence if necessary, in $L^1_{loc}$ and a.e.\\
Similarly to the case of BV data (\cite{Br00}), we can establish that the
sequence $(u^n)_{n\geq 0}$ is bounded in $ Lip^s([0,+\infty[_t,\,L^{1/s}_{loc}(\R_x,\R))$: this is
enough to get the convergence a.e. in $[0,+\infty[_t\times\R$ of some subsequence (still noted
$(u^n)$) towards a function $u$, entropy solution of the initial problem. Lastly Proposition
\ref{TVSCI}
ensures that for all $t\geq 0$ and $n\in\N$, $TV^su(t,\cdot)\leq TV^s u^n(t,\cdot)\leq TV^s u_0$,
thus Theorem \ref{stab} holds.
\cqfd


\section{Smoothing effect for nonlinear degenerate convex fluxes } \label{spsa}

First  we define the degeneracy of a nonlinear flux.  Then we obtain a
 smoothing effect in the spirit of P.-D. Lax \cite{La57} and O. Oleinik \cite{O57}. 
 Finally, we study the asymptotic behavior of entropy solutions as \cite{LP84}. 
There is two main tools: the Lax-Oleinik formula and the $BV^s$ spaces. We refer the reader
to the book of P.-D. Lax \cite{La06} for these results in the case of uniformly convex flux and also
to \cite{E98} for detailed proofs.

\subsection{Degenerate nonlinear flux}

\begin{definition}[degeneracy of a nonlinear  convex flux] \label{defdege}
\alali 
 Let $f$ belong to $C^1(I,\R)$ where $I$ is an interval of $\R$.  We say that the degeneracy of $f$ 
on $I$ is at least $p$ if the continuous derivative  $a(u)=f'(u)$ satisfies: 
\begin{equation}\label{defdeg}
   0  < \inf_{I \times I} \frac{|a(u) - a(v)|}{|u-v|^p}
\end{equation} 
We call the lowest real number $p$, if it exists, the degeneracy measurement of uniform convexity on
$I$.  If there is no $p$ such that  \eqref{defdeg} is satisfied, we set $p=+ \infty$.\\
Let $f\in C^2(I)$. We say that a real number $y\in I$ is a degeneracy point of $f$ in $I$ if
$f''(y)=0$ (i.e. $y$ is a critical point of $a$).
 
\end{definition}
If $f \in C^2(I)$ we can see easily that $p \geq 1$. 

\begin{remark}\label{strict}
Condition \eqref{defdeg} implies the strict monotonicity of $a(.)$ and then the strict 
convexity or concavity of the flux, but it is more general than the uniform convex case studied by
P.-D. Lax in \cite{La57}. 
Indeed, \eqref{defdeg} allows $f''$ to vanish as one can see below with the power law flux function.
\end{remark}
We give some examples  to illuminate this notion. 

\begin{ex}{\rm\textbf{Uniformly convex function:}  $ \inf f''> 0$.}\\
The degeneracy is $p=1$. 
\end{ex}
This is the basic example studied by P.-D. Lax \cite{La57} with $I=\R$. 
 
\begin{ex}{\rm\textbf{Linear flux.}} The degeneracy is  $p=+ \infty$. 
\end{ex}

\begin{ex}{\rm\textbf{Power convex functions  $f(u)=\ds \frac{ |u|^{1+\alpha}}{1+\alpha}$, $\alpha >
0$.}}\\
Let $I=[0,1]$, then $y=0$ is a degeneracy point and   the degeneracy of $f$ in $I$ is $p=
\max(1,\alpha)$.
\end{ex}
This example is the basic example to obtain all the  finite degeneracy  $p\geq 1$. 
\medskip

\pr 
the computation of $p$ is straightforward.  The case $\alpha < 1$ is left to the reader.
The case $\alpha=1$ corresponds to the Burgers flux, the simplest example of an uniformly strictly 
convex flux. Let us study the more interesting case $\alpha> 1$.
It is clear that $p \geq \alpha$, else the fraction of Inequality \eqref{defdeg} vanishes for $v=0$
and $u \rightarrow 0$. It suffices to study the case $p=\alpha$.
Let $R(u,v) = \ds  \frac{|u^\alpha - v^\alpha|}{|u-v|^\alpha}$  for $u \neq v$. 
It suffices to study the case $u < v$ by symmetry: with $v=u+h, \; h>0$, 
$R(u,u+h) = \ds \frac{ (u+h)^\alpha - u^\alpha}{h^\alpha}= \phi(y) = (y+1)^\alpha - y^\alpha $,
where $y =\displaystyle \frac{u}{h} \in [0,+\infty[$. Then $\displaystyle \inf_{y \geq  0} \phi(y)=
2^{1 -\alpha}>0$ which is enough to conclude. 
\cqfd

\begin{ex}{\rm\textbf{Smooth degenerate convex flux.}}
\alali
Let $K$ be a compact interval,  $f\in C^\infty(K,\R)$ and let $a=f'$ be an increasing function.
We define classically  the valuation of  $f''$   by:  $$ val[a](u)= \min \left \{k\geq 1, \frac{d^k
a}{d u^k}(u) \neq 0 \right\} \in \{1,2,\dots\} \cup \{+ \infty\}$$ then the degeneracy of $f$ on $K$
is $p = \ds \max_K val[a]$.
\end{ex}
We say that the flux is nonlinear if $p$ is finite.  

This general example has been studied recently for the multidimensional case in \cite{BJ10,Ju11}. 
These  examples allow to compute the parameter of degeneracy of  any smooth flux  given in the 
paper of P.-L. Lions, B. Perthame and E. Tadmor \cite{LPT94}.
\medskip

\pr 
In the one dimensional case, the computation is easier. We give a simple proof for a nonlinear flux,
i.e. the valuation is finite for each point of $K$.
Let $R(u,v) = \ds  \frac{|a(u) - a(v)|}{|u-v|^p}$ for $u \neq v$.  
Since $R$ is  a continuous function on $u \neq v$, positive outside the diagonal $\{u=v\}$, 
it suffices to study $R$ on the diagonal. 
Let  $k$ be $ val[a](u)$,   $ \ds R(u,u)= \left \{ \begin{array}{ccc}  0   &  if   &  k > p, \\ 
                                                                                                   
|a^{(k)}(u)| & if &  k=p, \\
+\infty   & if &  k  < p.    
\end{array}         
\right .$ 
\\
So the lowest $p$ in the neighborhood of $u$ is $p= val[a](u)$. 
Notice that the valuation is upper  semi-continuous. So the maximum of the valuation on the compact
$K$ exists and it is the lowest $p$  satisfying  Definition \ref{defdege} .
\cqfd

\subsection{Smoothing effect}  \label{sse}

We generalize the Oleinik one sided Lipschitz condition \cite{O57}  to define an entropy solution on
the scalar conservation law \eqref{eqscl} and we prove that the Lax-Oleinik formula yields such
condition for degenerate convex flux.

\begin{definition}[One sided H\"{o}lder condition]~\\
Let $f$ be a degenerate  convex flux.
Let $p \geq 1$ be a degeneracy pa\-ra\-meter of $f$ on an interval $I$,  and $0< s = \dfrac{1}{p} \leq 1$. 
Let be $u$ a  weak solution of \eqref{eqscl}.  Assume that  $u$ belongs to $I$. This solution is
called an entropy solution if  for some positive constant $c$, for all $t>0$ and  for almost all $(x,y)$ such that  $x <y$ we have
\begin{eqnarray} \label{os}
u(t,y)-u(t,x) &  \leq &  c  \, \frac{(y-x)^s}{t^s}.
\end{eqnarray}
\end{definition}
If $-f$ is convex then we replace in Inequality \eqref{os} $u$ by $-u$.
\\
As usual, the one sided condition implies the Lax entropy condition \cite{D00}.

\begin{theorem}[$BV^s$ smoothing effect  for  degenerate convex flux] \label{end}
\alali
Let $K $ be the compact interval $[-M,M]$.
Let  $u_0$  belong to $ L^\infty(\R)$,  $f \in C^1(\R,\R)$  and let $u$
be the unique entropy solution  on $]0,+\infty[_t\times\R_x$ of the scalar conservation law
\eqref{eqscl} satisfying the one sided condition \eqref{os}.
Let $p$ be a degeneracy parameter of $f$ on $K$  and  $0 < s = \ds \frac{1}{p} \leq 1$. 
\\
If  $p$ is finite and $|u_0| \leq M$  then  $u \in \mbox{\rm
Lip}^s(]0,+\infty[_t,L_{loc}^{1/s}(\R_x,\R))$ and $$\forall t>0, \, u(t,.) \in BV^s_{loc}(\R).$$
\\
If $u_0$ is compactly supported then $u \in \mbox{\rm Lip}^s(]0,+\infty[_t,L^{1/s}(\R_x,\R))$ and
there exists a constant $C$ such that $$TV^s u(t,.) \leq C \left( 1+ \displaystyle \frac{1}{t}
\right).$$
\end{theorem}

\begin{remark}
This entropy solution is the unique Kruzhkov entropy solution. It is well known for an uniformly
convex flux \cite{D00},
for the degenerate convex case see (\cite{Castel}.
\end{remark}

\begin{remark}
This theorem  gives the regularity conjectured by P.-L. Lions, B. Per\-thame and E. Tadmor in
\cite{LPT94} for a non linear convex flux. This conjecture was stated in Sobolev spaces.
The $W^{s,1}$ regularity with only $L^\infty $ initial data    was first proved in \cite{Ja09}.
We get  the best $W^{s,p}$ regularity. Indeed
by Proposition \ref{pSobolevBVs},   this $BV^s$ regularity gives a  $W^{s',1/s}$ smoothing effect for all
$s'<s$.
\end{remark}

\begin{remark}
We cannot expect a better regularity. Indeed, C. De Lellis and M.
Westdickenberg give in  \cite{DW} a piecewise smooth entropy solution which  does not belong to  
$W^{s,1/s}$. 
Recently, in \cite{Castel, CJ1},  another examples, with   continuous functions,  are built.
 Indeed for each  $\tau >s$ , there exists a smooth solution  which belongs to  $BV^s$  but not to  $BV^{\tau}$. 
\end{remark}

\begin{remark}
For  solutions with bounded entropy production and uniform convex flux,  the optimal smoothing effect is reached in \cite{GP}.
This class of solutions is larger than the class of entropy solutions. 
The optimal exponent is only $s=1/3$ (\cite{DW,GP}) instead of $s=1$ for uniformly convex fluxes. 

\end{remark}


\pr 
we first recall the  Lax-Oleinik formula for  a general convex flux without assuming the uniform
convexity. We assume only \eqref{defdeg}.  With such  an assumption the Lax-Oleinik formula is
still valid  (\cite{Castel}).
We know, thanks to Remark \ref{strict}, that the function $a$ (or $-a$), is  increasing. We assume
here that  the function  $a$ is   increasing on $K$.  We can easily extend $a$ continuously on $\R$
with the same degeneracy parameter $p$ (using a suitable translated power function) then the
function $a$ admits the inverse function $b$ on $\R$.
The entropy solution is then given for all $t$ and almost all $x$ by  the  Lax-Oleinik formula:  
\begin{equation}\label{LOF}
\ds  u(t,x)= b \left ( \frac{x-y}{t}\right)
\end{equation} 
where $y=y(t,x)$ minimizes, for $t$ and $x$ fixed, the function  
$$\ds G(t,x,y)=U_0(y) + t \, h\left(\frac{x-y}{t}\right)$$
 with $\ds U_0(y)= \int_0^y u_0(x)\,dx$,  $a(0)=c$,  $ \ds h(u) = \int_c^u b(v)\,dv$. 

Geometrically, $y(t,x)$ has a simple interpretation. The function $u(.,.)$ is constant on the
cha\-racteristic $x=y + t a(u_0(y))$: $u(t,x)=u_0(y)$ (before the formation of a shock).
Indeed $a(u_0(y))=\ds \frac{x-y}{t}$, so  $b(a(u_0(y))) = b\left( \ds\frac{x-y}{t} \right)=
u_0(y)=u(t,x) $. The key point of the formula \eqref{LOF}  is that $y(t,x)$ minimizes an explicit  
function, namely $y \mapsto G(t,x,y)$. Consequently $y(t,x)$ is not so far from $x$, more precisely:
  \begin{equation}\label{yx}
 | x- y(t,x)|  \leq  t  \, \sup_K |a|. 
\end{equation} 
Moreover, if $x_1 < x_2$ then $y(t,x_1) \leq y(t,x_2)$,  (\cite{La06,E98,Castel}).

Condition (\ref{defdeg})  implies that $b$ belongs to  $C^s(\R,\R)$   with $s = 1/p $. Indeed
we have  for all $U,V \in a(K)$,  with $u=b(U)$ and $v=b(V)$:
$$ \ds  \frac{|b(U) - b(V)|}{|U-V|^s }=  \frac{|u - v|}{|a(u)-a(v)|^s } = \left( \frac{|u -
v|^p}{|a(u)-a(v)| } \right)^s \leq  D^s= \frac{1}{C^s}$$
where 
$ \ds 0 < C =    \inf_{K \times K} \frac{|a(u) - a(v)|}{|u-v|^p} $.

We are now able to prove the $BV^s$ smoothing effect.  Fix $T>0$  and $I=[a,b]$: we want to bound
$TV^s u \{I\}$. Let $x_1,\,x_2 \in I$ and $y_i = y(t,x_i)$, then
\begin{eqnarray*}
 |u(T,x_1)-u(T,x_2)|^p &=& \ds  \left |  b \left ( \frac{x_1-y_1}{T}\right) -  b \left (
\frac{x_2-y_2}{T}\right) \right|^p\\
& \leq & \left(  D^s  \left |   \frac{x_1-y_1}{T} -\frac{x_2-y_2}{T} \right|^s\right)^p.
\end{eqnarray*}
The condition  $s p=1$ yields 
$$ |u(T,x_1)-u(T,x_2)|^p \leq  \ds D  \left |   \frac{x_1-x_2}{T} \right| + D  \left |
\frac{y_1-y_2}{T} \right|.$$
We  now compute  $TV^s u \{\sigma\}$ for a subdivision $\sigma =\{x_0 < x_1< \cdots < x_n\}$ of $I$.
Then
\begin{equation} \label{TVsLOF}
TV^s u \{\sigma\} \leq \frac{D}{T} (x_n -x_0 + y_n-y_0) \leq  \frac{D}{T} (2(b-a)  + T  \sup_K|a|).
\end{equation}
Then $TV^s u \{I\}$  keeps the same bound.
\\
We can precise the previous bounds.
First, we obtain  the one sided H\"{o}lder condition \eqref{os}, which implies that the solution is
an
entropy solution. We know that  if $x_1 < x_2$ then $$y_1=y(t,x_1) \leq y(t,x_2)=y_2.$$
Moreover
$$u(t,x_2)-u(t,x_1) =b \left( \frac{x_2-y_2}{t} \right)-b \left( \frac{x_1-y_1}{t} \right) \leq 
b \left( \frac{x_2-y_2}{t} \right)-b \left( \frac{x_1-y_1}{t} \right)$$ because $b$ is increasing.
But,  $b(\frac{x_2-y_2}{t})-b(\frac{x_1-y_1}{t}) \geq 0$ because $x_2 \geq x_1$. Then
$$u(t,x_2)-u(t,x_1)  \leq \left |   b \left( \frac{x_2-y_2}{t} \right) -b \left( \frac{x_1-y_1}{t}
\right) \right|
\leq D^s \left |   \frac{x_2-x_1}{t} \right| ^s = D^s  \frac{(x_2-x_1)^s}{t^s}.$$
We can  improve the $TV^s$ bound for a compactly supported initial data.
For any $t$, the solution stays compactly supported (but the size of this support depends on  $t$).
Fix $T=1$. Inequality \eqref{TVsLOF}  gives $u \in BV^s(\R)$.
\\
For $t \geq T$, Theorem \ref{stab} implies $TV^su(t,.) \leq TV^su(T,.)=C_1$ and 
for $0<t \leq T$, Inequality \eqref{TVsLOF} implies $TV^s u(t,.) \leq C_0(1+\ds\frac{1}{t})$,
then $\forall t>0$, $TV^s u(t,.) \leq C(1+\ds\frac{1}{t})$.
\\
Theorem \ref{stab}  shows that $u \in \mbox{\rm Lip}^s(]0,+\infty[_t,L^{1/s}(\R_x,\R))$.
\\
Fo the general case, the estimate is only locally valid with respect to the space variable. \cqfd

\begin{proposition}\label{pend}
The unique entropy solution  of  Theorem \ref{end} satisfies the folowing decay
$ \ds TV^s_+ u(T,.) \{[a,b]\} \leq   D\frac{|b-a|}{T}$ for some positive constant $D$.
\end{proposition}

\pr 
it is a direct consequence of the one sided condition \eqref{os}.
\cqfd


\subsection{Asymptotic behavior of entropy solutions}

The smoothing effect is sometimes  related  to the asymptotic behavior for large time (\cite{La57,La74,La06}).
We investigate briefly classical decays under assumption \eqref{defdeg}.
Indeed the decay of the solution with compact support   depends on one more  parameter.

\begin{theorem}[Decay for large time] \label{thdecay}
\alali
Let be  $u_0 \in L^1 \cap L^\infty(\R)$, $|u_0|\leq M$, $K=[-M,M]$.
Assume that  $f \in C^1(\R,\R)$ satisfies Condition \eqref{defdeg} with $a=f'$, $p$ the degeneracy
of $f$ on $K$, $s=1/p$.
Let  $u$ be  the unique entropy solution  on $]0,+\infty[_t\times\R_x$ of \eqref{eqscl} and   $b$
the inverse function of the function $a$ on $a(K)$. 
\\
If  there exists  $q>0$  such that 
\begin{equation} \label{condb}
0  < \inf_{U \in a(K)}  \frac{|b(U-a(0))|}{|U|^q} 
\end{equation}
then  there exists $C>0$ such that 
$$
 |u(t,x)|   \leq  \frac{C}{t^{d}}, \qquad    d =\frac{s}{1+q} .
$$
\end{theorem}
Originally, P.-D. Lax found this optimal decay in the 50' for  strictly convex flux  with  $d
=\ds\frac{1}{2}$ since   $s=q=1$ \cite{La57}. 
\\
For power function $f(u)= |u|^{1+\alpha}$, $\alpha > 1$,  we have  $d =\ds\frac{s}{1+s} <
\ds\frac{1}{2}$ since   $s=\ds q=1/\alpha$.
For the simplest degenerate convex case:
the cubic convex flux, we only have $d=\ds\frac{1}{3}$. This decay is slower than classical Lax
decay which is $\ds\frac{1}{\sqrt{t}}$.
\\
Remark that $q \geq s$.
Assume that without  loss of generality $a(0)=0$ and $J=a(K)=[0,1]$.
Then $| b(x)-b(y) | \leq C | x-y|^s$ since $b$ is in $C^s(J)$.
Moreover,  there exists $D>0$ such that $D|x|^q \leq |b(x)|$ by \eqref{condb},
so $Dx^q \leq C x^s$ on $[0,1]$ then we have $q \geq s$.
\\
We give some examples with $q>s$.
On $K=J=[0,1]$, with $f(u)=\ds\frac{u^{1+\alpha}}{1+\alpha}, \, 0 < \alpha < 1$ we have $s=1,\,
q=\ds\frac{1}{\alpha}>1=s$.
\medskip

\pr 
the proof is a slight modification of the original  Lax's proof, \cite{La06}. 
We use the Lax-Oleinik formula with the notations of the proof of Theorem \ref{end},
so we have to extend the function $a$ on $\R$.
We have

\begin{eqnarray*}
\forall y \in \R, \, -d_2 &\leq&  U_0(y)=\int_0^y u_0(x) \, dx\\
&\leq& d_2 = \max \left
( \int_0^{+ \infty} |u_0(x)|dx, \int^0_{- \infty} |u_0(x)|dx \right). 
\end{eqnarray*}

Notice that $\min G \leq d_2$. Since
$h$ is a convex nonnegative function which vanishes only at $ c= a(0)$, it suffices to take  $y=x-c
\, t$ so $G(t,x,y)=U_0(y)$. Integrating  Inequality \eqref{condb}, there exists a constant 
$d_1>0$ such that  for $z \in J=a(K)$, $$h(z) \geq d_1 | z -c|^{1+q}. $$
\\
Let $y=y(t,x)$ be the minimizer of $G(t,x,.)$.
\\
Notice that $ \ds\frac{x-y}{t}  \in J$ since $ x -y =t  \, a(u_0(y))$.
Now, we have   the inequality
$$\ds d_2 \geq  G(t,x,y) \geq -d_2 + t d_1\left | \frac{x-y}{t} - c \right |^{1+q}$$ then
$$ \ds  2 \, d_2 \geq   t \, d_1\left | \frac{x-y}{t} - c \right |^{1+q}$$ 
and then 
\begin{eqnarray} \label{ineqhy}
\ds  \left (\frac{2 d_2}{  d_1 t } \right)^{1/(1+q)} & \geq &  \left | \frac{x-y}{t} - c \right | .
\end{eqnarray}
Since $b \in C^s$, we have $ |b(z)| = |b(z) - b(c)| \leq D^s |z-c|^s$ .
\\
The Lax-Oleinik formula \eqref{LOF}   and  Inequality \eqref{ineqhy} conclude the proof: 
\begin{eqnarray*}
\ds  |u(t,x)| =\left | b \left ( \frac{x-y}{t} \right)  \right |  &   \leq &  \ds    D^s  \left
(\frac{2 \, d_2}{  d_1 t } \right)^{s/(1+q)}.
\end{eqnarray*}
\cqfd

The periodic case is much simpler  and  only depends on the degeneracy of $f$.

\begin{theorem}[Decay for periodic solutions] \label{thdecayp}
\alali
Let $u_0$ be a P-periodic bounded function,   $m =\ds \frac{1}{P}\int_0^P u_0(x) dx $,  $|u_0|\leq
M$,\\ $K=[-M,M]$, let $u$ be the unique entropy solution  on $]0,+\infty[_t\times\R_x$ of
\eqref{eqscl}, $f \in C^1(\R,\R)$. If   the degeneracy $p$ of $f$ on $K$ is finite then there exists
a constant $C$ such that 
$$
 |u(t,x)- m |    \leq  \frac{C}{t^{s}}, \, s=\frac{1}{p}.
$$
\end{theorem}
For  uniform  convex flux we have   the classical case with   $s=1$, \cite{La57}. 
\\
For power function $f(u)= |u|^{1+\alpha}$ with $\alpha > 1$ we have  $s=\ds 1/\alpha$.
For instance, for the cubic convex flux, $s=\ds \frac{1}{2}$.   
\medskip

\pr
first notice that $u(t,.)$ is periodic with the same period $P$ and the same mean value $m$.
We have thanks to the one side condition \eqref{os} the inequality
$$u(t,y)-u(t,x) \leq C \frac{(y-x)^s}{t^s} \leq C \frac{P^s}{t^s}$$ for $0 \leq y-x \leq P$.
Assume that $m=0$ without  loss of generality.
Fix $x$. If $u(t,x)<0$, there exists $y$ in $[x,x+P]$ such that $u(t,y)>0$ since $m=0$.
Then 
$$|u(t,x)| \leq |u(t,x)| + |u(t,y)| \leq |u(t,y) - u(t,x)|= u(t,y) - u(t,x) \leq C
\frac{P^s}{t^s}.$$
The same argument holds if $u(t,x)>0$, which concludes the proof.
\cqfd



\end{document}